\documentclass[12pt,a4paper]{amsart}

\usepackage{amssymb}
\usepackage{amsmath}
\usepackage{amsthm}
\usepackage{xypic}

\theoremstyle{plain}
\newtheorem{theorem}{Theorem}
\newtheorem{proposition}[theorem]{Proposition}

\theoremstyle{remark}
\newtheorem{remark}[theorem]{Remark}
\theoremstyle{definition}
\newtheorem{definition}[theorem]{Definition}

\renewcommand{\baselinestretch}{1.2}

\def\varinjlim_#1{\lim\limits_{\longrightarrow\atop{#1}}}

\def\End{\mathop{\rm End}\nolimits}

\def\Aut{\mathop{\rm Aut}\nolimits}

\def\Mor{\mathop{\rm Mor}\nolimits}
\def\Hom{\mathop{\rm Hom}\nolimits}
\def\id{\mathop{\rm id}\nolimits}

\def\U{\mathop{\rm U}\nolimits}
\def\PU{\mathop{\rm PU}\nolimits}

\def\B{\mathop{\rm B}\nolimits}

\def\K{\mathop{\rm K}\nolimits}
\def\diag{\mathop{\rm diag}\nolimits}
\def\pt{\mathop{\rm pt}\nolimits}
\def\tr{\mathop{\rm tr}\nolimits}

\def\im{\mathop{\rm im}\nolimits}

\def\Hom{\mathop{\rm Hom}\nolimits}

\def\Gr{\mathop{\rm Gr}\nolimits}
\def\Fr{\mathop{\rm Fr}\nolimits}

\def\BU{\mathop{\rm BU}\nolimits}

\def\BPU{\mathop{\rm BPU}\nolimits}

\def\EG{\mathop{\rm EG}\nolimits}

\def\SU{\mathop{\rm SU}\nolimits}
\def\EPU{\mathop{\rm EPU}\nolimits}

\def\BBU{\mathop{\rm BBU}\nolimits}
\def\BBSU{\mathop{\rm BBSU}\nolimits}
\def\BSU{\mathop{\rm BSU}\nolimits}

\def\Fred{\mathop{\rm Fred}\nolimits}

\def\Id{\mathop{\rm Id}\nolimits}

\begin{document}

\title{On $K$-theory automorphisms related to bundles of finite order}
\author{A.V. Ershov}

\email{ershov.andrei@gmail.com}

\begin{abstract}
In the present paper we describe
the action of (not necessarily line) bundles
of finite order
on the $K$-functor in terms of classifying spaces.
This description might provide with an approach for
more general twistings in $K$-theory than ones related
to the action of the Picard group.
\end{abstract}

\date{}
\maketitle {\renewcommand{\baselinestretch}{1.0}


\section*{Introduction}

The complex $K$-theory is a generalized cohomology theory
represented by the $\Omega$-spectrum $\{ K_n\}_{n\geq 0}$, where
$K_n=\mathbb{Z}\times \BU$, if $n$ is even and $K_n={\rm U}$ if
$n$ is odd. $K_0=\mathbb{Z}\times \BU$ is an $E_\infty$-ring
space, and the corresponding space of units $K_\otimes$ (which is an infinite loop space)
is $\mathbb{Z}/2\mathbb{Z}\times
\BU_\otimes$, where $\BU_\otimes$ denotes the space $\BU$ with the
$H$-space structure induced by the tensor product of virtual bundles of virtual dimension $1$. Twistings of the
$K$-theory over a compact space $X$ are classified by homotopy classes of maps
$X\rightarrow {\rm
B}(\mathbb{Z}/2\mathbb{Z}\times \BU_\otimes)\simeq
\K(\mathbb{Z}/2\mathbb{Z},\, 1)\times \BBU_\otimes$ (where $\rm B$
denotes the functor of classifying space).
The theorem that
$\BU_\otimes$ is an infinite loop space was proved by G. Segal \cite{Segal}. Moreover, the spectrum
$\BU_\otimes$ can be decomposed as follows:
$\BU_\otimes = \K(\mathbb{Z},\, 2)\times
\BSU_\otimes$.
This implies that the twistings in the
$K$-theory can be classified by homotopy classes of maps
$X\rightarrow \K(\mathbb{Z}/2\mathbb{Z},\, 1)\times
\K(\mathbb{Z},\, 3)\times \BBSU_\otimes$.
In other words, for a compact space $X$ the twistings correspond to elements in
$H^1(X,\, \mathbb{Z}/2\mathbb{Z})\times H^3(X,\, \mathbb{Z})\times
[X,\, \BBSU_\otimes],\; [X,\, \BBSU_\otimes]= bsu^1_\otimes(X),$
where $\{ bsu^n_\otimes\}_n$ is the generalized cohomology theory
corresponding to the infinite loop space $\BSU_\otimes$.

The twisted $K$-theory corresponding to the twistings coming from $H^1(X,\, \mathbb{Z}/2\mathbb{Z})\times H^3(X,\, \mathbb{Z})$ has
been intensively studied during the last decade, but not the general case (as far as the author knows).
It seems that the reason is that there is no known appropriate geometric model for ``nonabelian'' twistings from $[X,\, \BBSU_\otimes].$
In this paper we make an attempt to give such a model for elements of finite order in $[X,\, \BBSU_\otimes].$
In particular, we are based on the model of the $H$-space $\BSU_\otimes$ given by the infinite matrix grassmannian $\Gr$ \cite{Ers1}
(see also subsection 4.5 below).

A brief outline of this paper is as follows. In section 1 we recall the well-known result that
the action of the projective unitary group of the separable Hilbert space $\PU({\mathcal H})$
on the space of Fredholm operators $\Fred({\mathcal H})$ (which is the
representing space of $K$-theory) by conjugation corresponds
to the action of the Picard group $Pic(X)$ on $K(X)$ by group automorphisms (Theorem \ref{lineth}).
The key result of section 2 is Theorem \ref{mth} which is in some sense a counterpart of Theorem
\ref{lineth}. Roughly speaking, it asserts that in terms of representing space $\Fred({\mathcal H})$ the tensor multiplication of
$K$-functor by (not necessarily line) bundles of finite order $k$ can be described by some maps
$\gamma_{kl^m,\, l^{n-m}}^\prime \colon \Fr_{kl^m,\, l^{n-m}}\times \Fred({\mathcal H}) \rightarrow \Fred({\mathcal H}),$
where $\Fr_{kl^m,\, l^{n-m}}$ are the spaces parametrizing unital $*$-homomorphisms of matrix algebras $M_{kl^m}(\mathbb{C})\rightarrow M_{kl^n}(\mathbb{C})$. Then by arranging these maps
$\gamma_{kl^m,\, l^{n-m}}^\prime$ we should construct an action of the $H$-space $\varinjlim_n \Fr_{kl^n,\, l^n}$ on
$\Fred({\mathcal H})$. We also have shown that $\varinjlim_n \Fr_{kl^n,\, l^n}$ is a well-pointed grouplike topological monoid
and therefore there exists the classifying space $\B \varinjlim_n \Fr_{kl^n,\, l^n}$ (see subsection 3.2).
More precisely, we consider the direct limit of matrix algebras
$M_{kl^\infty}(\mathbb{C})=\varinjlim_mM_{kl^m}(\mathbb{C})$ and the monoid of its unital endomorphisms.
We fix an increasing filtration $A_{kl^m}\subset A_{kl^{m+1}}\subset \ldots,\quad
A_{kl^m}=M_{kl^m}(\mathbb{C})$ in
$M_{kl^\infty}(\mathbb{C})$ such that $A_{kl^{m+1}}=M_l( A_{kl^m})$.
Then we consider endomorphisms of $M_{kl^\infty}(\mathbb{C})$ that are induced by
unital homomorphisms of the form $h\colon A_{kl^m}\rightarrow A_{kl^n}$ (for some $m,\, n$), i.e. have the form
$M_{l^\infty}(h).$
Such endomorphisms form the above mentioned topological monoid which is
homotopy equivalent to the direct limit
$\Fr_{kl^\infty,\, l^\infty}:=\varinjlim_{m,\, n}\Fr_{kl^m,\, l^n}$,
which is not contractible provided $(k,\, l)=1$.
Moreover, the automorphism subgroup in the monoid (corresponding to $n=0$) is
$\varinjlim_m\PU(kl^m).$ Furthermore, this monoid naturally acts on the space of Fredholm operators
and this action corresponds to the tensor multiplication of the $K$-functor by bundles
of order $k$. In subsection 3.3 we also sketch the idea of the definition of the corresponding version of the twisted $K$-theory.

Roughly speaking, the ``usual'' (Abelian) twistings of order $k$ correspond to
the group of automorphisms while the nonabelian ones correspond to the monoid of endomorphisms of $M_{kl^\infty}(\mathbb{C})$.
Note that these endomorphisms act on the localization of the space of Fredholm operators over $l$ by homotopy auto-equivalences,
i.e. they are invertible in the sense of homotopy.

Although some technical difficulties remain we hope that this approach will be useful in order to
define a general version of the twisted $K$-theory.

\section{$K$-theory automorphisms related to line bundles}

In this section we describe well-known results about the action of
$Pic(X)$ on the group $K(X)$. We also consider the special case
of the subgroup of line bundles of finite order.

Let $X$ be a compact space,  $Pic(X)$ its Picard group consisting of isomorphism classes of line bundles
with respect to the tensor product.
The Picard group is represented by the $H$-space $\BU(1)\cong \mathbb{C}P^\infty \cong \K(\mathbb{Z},\, 2)$
whose multiplication is given by the tensor product of line bundles or (in the appearance of the Eilenberg-MacLane space) by addition
of two-dimensional integer cohomology classes.
In particular, the first Chern class $c_1$ defines the isomorphism
$c_1 \colon Pic(X)\stackrel{\cong}{\rightarrow}H^2(X,\, \mathbb{Z}).$
The group $Pic(X)$ is a subgroup of the multiplicative group of the ring $K(X)$
and therefore it acts on $K(X)$ by group automorphisms.
This action is functorial on $X$ and therefore it can be described in terms of classifying spaces
(see Theorem \ref{lineth}).

As a representing space for the
$K$-theory we take $\Fred({\mathcal H}),$ the space of Fredholm operators in the separable Hilbert space
${\mathcal H}$. It is known
\cite{AS} that for a compact space $X$ the action of
$Pic(X)$ on $K(X)$ is induced by the conjugate action
$$
\gamma \colon \PU({\mathcal H})\times \Fred({\mathcal
H})\rightarrow \Fred({\mathcal H}),\; \gamma(g,\, T)=gTg^{-1}
$$
of $\PU({\mathcal H})$ on $\Fred({\mathcal H}).$ More precisely, there is the following theorem
(recall that $\PU({\mathcal H})\simeq \mathbb{C}P^\infty \simeq
\K(\mathbb{Z},\, 2)$).

\begin{theorem}
\label{lineth}
If $f_\xi \colon X\rightarrow
\Fred({\mathcal H})$ and $\varphi_\zeta \colon X\rightarrow
\PU({\mathcal H})$ represent $\xi \in K(X)$ and $\zeta \in Pic(X)$
respectively, then the composite map
\begin{equation}
\label{compmapconj}
X\stackrel{\diag}{\longrightarrow}X\times
X\stackrel{\varphi_\zeta \times
f_\xi}{\longrightarrow}\PU({\mathcal H})\times \Fred({\mathcal H})
\stackrel{\gamma}{\rightarrow}\Fred({\mathcal H})
\end{equation}
represents $\zeta \otimes \xi \in K(X)$.
\end{theorem}
{\noindent \it Proof \;} see \cite{AS}.$\quad \square$

\smallskip

It is essential for the theorem that the group $\PU({\mathcal H})$, on the one hand having the homotopy type of $\mathbb{C}P^\infty$ is the base
of the universal $\U(1)$-bundle (which is related to the exact sequence of groups $\U(1)\rightarrow \U({\mathcal H})\rightarrow \PU({\mathcal H})$,
because $\U({\mathcal H})$ is contructible in the considered norm topology), on the other hand being a group acts in the appropriate way
on the representing space of $K$-theory (the space of Fredholm operators).

Then in order to define the corresponding version of the twisted $K$-theory
one considers the $\Fred({\mathcal H})$-bundle
$\widetilde{\Fred}({\mathcal H})\rightarrow \BPU({\mathcal H})$ associated (by means of the action $\gamma$) with
the universal $\PU({\mathcal H})$-bundle over the classifying space
$\BPU({\mathcal H})\simeq \K(\mathbb{Z},\, 3)$, i.e. the bundle
\begin{equation}
\label{bunndd}
\diagram
\Fred({\mathcal H}) \rto & \EPU({\mathcal H}){\mathop{\times}\limits_{\PU({\mathcal H})}}\Fred({\mathcal H}) \dto \\
& \BPU({\mathcal H}). \\
\enddiagram
\end{equation}

Then for any map $f\colon X\rightarrow \BPU({\mathcal H})$ the corresponding twisted $K$-theory $K_f(X)$
is the set (in fact the group) of homotopy classes of sections
$[X,\, f^*\widetilde{\Fred}({\mathcal H})]^\prime$ of the pullback bundle (here $[\ldots ,\, \ldots ]^\prime$
denotes the set of fiberwise homotopy classes of sections).
The group $K_f(X)$ depends up to isomorphism only on the homotopy class $[f]$ of the map $f$,
i.e. in fact on the corresponding third integer cohomology class.

In this paper we are interested in the case of bundles
(more precisely, of elements in $bsu_\otimes^0$) of finite order, therefore let us consider separately
the specialization of the mentioned result to the case of line bundles of order
$k$ in $Pic(X).$
For this we should consider subgroups $\PU(k)\subset\PU({\mathcal H})$.
Let us describe the corresponding embedding.

Let ${\mathcal B}({\mathcal H})$ be the algebra of bounded operators on the separable Hilbert space
${\mathcal H}$, $M_k({\mathcal B}({\mathcal
H})):=M_k(\mathbb{C}){\mathop{\otimes}\limits_{\mathbb{C}}}{\mathcal
B}({\mathcal H})$ the matrix algebra over ${\mathcal B}({\mathcal H})$
(of course, it is isomorphic to ${\mathcal B}({\mathcal H})$).
Let $\U_k({\mathcal H})\subset M_k({\mathcal B}({\mathcal H}))$
be the corresponding unitary group (which is isomorphic to $\U({\mathcal H})$). It acts on $M_k({\mathcal B}({\mathcal
H}))$ by conjugations (which are $*$-algebra isomorphisms), moreover, the kernel of the action is the center, i.e. the subgroup
of scalar matrices $\cong \U(1).$ The corresponding quotient group we denote by
$\PU_k({\mathcal H})$ (of course, it is isomorphic to $\PU({\mathcal H})$).

$M_k(\mathbb{C})\otimes \Id_{{\mathcal B}({\mathcal H})}$ is a {\it $k$-subalgebra} (i.e. a unital $*$-subalgebra isomorphic to $M_k(\mathbb{C})$) in
$M_k({\mathcal B}({\mathcal H}))$. Then $\PU(k)\subset \PU_k({\mathcal H})$
is the subgroup of automorphisms of this $k$-subalgebra.
Thereby we have defined the injective group homomorphism
$$
i_k\colon \PU(k)\hookrightarrow \PU_k({\mathcal H})
$$
induced by the group homomorphism $\U(k)\hookrightarrow \U_k({\mathcal H}),\: g\mapsto
g\otimes \Id_{{\mathcal B}({\mathcal H})}$.

Let $[k]$ be the trivial $\mathbb{C}^k$-bundle over $X$.

\begin{proposition}
\label{ktors}
For a line bundle
$\zeta \rightarrow X$ satisfying the condition
\begin{equation}
\label{condtriv}
\zeta \otimes [k]=\zeta^{\oplus k}\cong X\times \mathbb{C}^k
\end{equation}
the classifying map $\varphi_{\zeta} \colon
X\rightarrow \PU_k({\mathcal H})\cong \PU({\mathcal H})$ can be lifted to a map
$\widetilde{\varphi}_\zeta \colon X\rightarrow \PU(k)$ such that
$i_k\circ \widetilde{\varphi}_\zeta \simeq \varphi_{\zeta}$.
\end{proposition}
{\noindent \it Proof.}
Consider the exact sequence of groups
\begin{equation}
\label{flb}
1\rightarrow \U(1)\rightarrow
\U(k)\stackrel{\chi_k}{\rightarrow}\PU(k) \rightarrow 1
\end{equation}
and the fibration
\begin{equation}
\label{flb2}
\PU(k)\stackrel{\psi_k}{\rightarrow}\BU(1)\stackrel{\omega_k}{\rightarrow}\BU(k)
\end{equation}
obtained by its extension to the right.
In particular, $\psi_k \colon \PU(k)\rightarrow \BU(1)\simeq
\mathbb{C}P^\infty$ is the classifying map for the $\U(1)$-bundle $\chi_k$
(\ref{flb}).
It is easy to see that the diagram
$$
\diagram
\PU(k)\drto_{i_k}\rto^{\psi_k} & \BU(1)\\
& \PU_k({\mathcal H})\uto_{\simeq} \\
\enddiagram
$$
commutes.

Let $\zeta \rightarrow X$ be a line bundle satisfying the condition (\ref{condtriv}),
$\varphi_\zeta \colon X\rightarrow
\BU(1)$ its classifying map.
Since
$\omega_k$ (see (\ref{flb2})) is induced by taking the direct sum of a line bundle with itself $k$ times (and the extension of the structural group to $\U(k)$),
we see that $\omega_k\circ \varphi_\zeta \simeq *$. Now it is easy to see from exactness of (\ref{flb2}) that
$\varphi_\zeta \colon X\rightarrow
\BU(1)$ can be lifted to $\widetilde{\varphi}_\zeta
\colon X\rightarrow \PU(k). \quad \square$

\smallskip

Note that the choice of a lift
$\widetilde{\varphi}_\zeta$ corresponds to the choice of a trivialization
(\ref{condtriv}): two lifts differ up to a map
$X\rightarrow \U(k)$. Thus, a lift is defined up to the action of
$[X,\, \U(k)]$ on $[X,\, \PU(k)].$
The subgroup in $Pic(X)$ consisting of line bundles satisfying the condition (\ref{condtriv}) is
$\im \{ \psi_{k*}\colon [X,\, \PU(k)]\rightarrow [X,\,
\mathbb{C}P^\infty]\}$ or the quotient $[X,\, \PU(k)]/[X,\, \U(k)].$

Let $\Fred_k({\mathcal H})$ be the subspace of Fredholm operators in
$M_k({\mathcal B}({\mathcal H}))$. Clearly, $\Fred_k({\mathcal H})\cong \Fred({\mathcal H})$.
Acting on $M_k(\mathbb{C})$ by $*$-automorphisms, the group $\PU(k)$ acts on $M_k(\mathbb{C})\otimes \mathcal{B}(\mathcal{H})=M_k({\mathcal B}({\mathcal H}))$
through the first tensor factor.
Let $\gamma^\prime_k \colon \PU(k)\times \Fred_k({\mathcal H})\rightarrow \Fred_k({\mathcal H})$ be the restriction of this action on $\Fred_k({\mathcal H})$.
Then the diagram
$$
\diagram
\PU({\mathcal H})\times \Fred_k({\mathcal H}) \rto^{\qquad \gamma} & \Fred_k({\mathcal H}) \\
\PU(k)\times \Fred_k({\mathcal H}) \uto^{i_k\times \id} \urto_{\gamma^\prime_k} & \\
\enddiagram
$$
commutes. Now one can consider the $\Fred_k({\mathcal H})$-bundle
\begin{equation}
\label{linfin}
\begin{array}{c}
\diagram
\Fred_k({\mathcal H})\rto & \EPU(k){\mathop{\times}\limits_{\PU(k)}}\Fred_k({\mathcal H}) \dto \\
& \BPU(k) \\
\enddiagram
\end{array}
\end{equation}
associated by means of the action $\gamma^\prime_k$. This bundle is the pullback of (\ref{bunndd})
by $\B i_k$.

It is easy to see from the definition of the embedding
$i_k$ that the action $\gamma^\prime_k$ is trivial on elements of the form $k\xi.$
Indeed, a classifying map for
$k\xi$ can be decomposed into the composite
$X\stackrel{f_\xi}{\rightarrow}\Fred({\mathcal H})\stackrel{\diag}{\rightarrow}\Fred_k({\mathcal H}).$
From the other hand, $(1+(\zeta -1))\cdot k\xi =k\xi+0=k\xi$ or $\zeta \otimes ([k]\otimes \xi) =(\zeta \otimes [k]) \otimes \xi =[k]\otimes \xi.$

\begin{remark}
\label{projgract}
Note that if we choose an isomorphism ${\mathcal B}({\mathcal H})\cong M_{k^\infty}({\mathcal B}({\mathcal H}))$
and hence the isomorphism $\Fred({\mathcal H})\cong \Fred_{k^\infty}({\mathcal H}),$
we can define the limit action $\gamma^\prime_{k^\infty}\colon \PU(k^\infty)\times \Fred_{k^\infty}({\mathcal H})\rightarrow \Fred_{k^\infty}({\mathcal H}),$ etc.
\end{remark}

\section{The case of bundles of dimension $\geq 1$}

As was pointed out in the previous section, the group
$\PU({\mathcal H})$, from one hand acts on the representing space of $K$-theory
$\Fred({\mathcal H})$, from the other hand it is the base of the universal line bundle.
This two facts lead to the result that the action
of $\PU({\mathcal H})$ on $K(X)$ corresponds to the tensor product by elements of the Picard group
$Pic(X)$ (i.e. classes of line bundles). This action can be restricted to subgroups
$\PU(k)\subset \PU({\mathcal H})$ which classify elements of finite order $k,\, k\in \mathbb{N}.$

In what follows the role of
groups $\PU(k)$ will play some spaces $\Fr_{k,\, l}$ (defined below).
From one hand, they ``act'' on $K$-theory (more precisely, their direct limit (which has the natural structure of an $H$-space)
acts), from the other hand, they are bases of some nontrivial
$l$-dimensional bundles of order $k$. We will show that their ``action'' on $K(X)$ corresponds to the tensor product by
those $l$-dimensional bundles (see Theorem \ref{mth}). The key result of this section is Theorem
\ref{mth} which can be regarded as a counterpart of Theorem \ref{lineth}.

Fix a pair of positive integers $k,\, l>1.$
Let ${\rm Hom}_{alg}(M_k(\mathbb{C}),\, M_{kl}(\mathbb{C}))$ be the space of all
unital $*$-homomorphisms $M_k(\mathbb{C})\rightarrow M_{kl}(\mathbb{C})$ \cite{Ers2}.
It follows from Noether-Skolem's theorem \cite{Pirce} that it can be represented in the form of a homogeneous space of the group $\PU(kl)$ as follows:
\begin{equation}
\label{reprfr}
{\rm Hom}_{alg}(M_k(\mathbb{C}),\, M_{kl}(\mathbb{C}))\cong \PU(kl)/(E_k\otimes \PU(l))
\end{equation}
(here $\otimes$ denotes the Kronecker product of matrices).
This space we shall denote by $\Fr_{k,\, l}$. We will be interested in the case $(k,\, l)=1$
(we have to impose a condition of such a kind to make the direct limit of the above spaces
noncontractible and the construction below nontrivial \cite{Ers1}).

\begin{proposition}
\label{interpret}
A map $X\rightarrow \Fr_{k,\, l}$ is the same thing as an embedding
\begin{equation}
\label{muu}
\mu \colon X\times M_k(\mathbb{C})\hookrightarrow
X\times M_{kl}(\mathbb{C}),
\end{equation}
whose restriction to a fiber is a unital $*$-homomorphism of matrix algebras.
\end{proposition}
{\noindent \it Proof \; } follows directly from the natural bijection
$$
\Mor (X\times M_k(\mathbb{C}),\, M_{kl}(\mathbb{C}))\rightarrow \Mor (X,\, \Mor (M_k(\mathbb{C}),\, M_{kl}(\mathbb{C}))).\quad \square
$$

\smallskip

Let $\Gr_{k,\, l}$ be the ``matrix grassmannian'', i.e. the space whose points parameterize
unital $*$-subalgebras in $M_{kl}(\mathbb{C})$ isomorphic to $M_{k}(\mathbb{C})$ (``$k$-{\it subalgebras}'') \cite{Ers1}.
It follows from Noether-Skolem's theorem \cite{Pirce} that
\begin{equation}
\label{reprmgr}
\Gr_{k,\, l}=\PU(kl)/(\PU(k)\otimes \PU(l)).
\end{equation}

The matrix grassmannian $\Gr_{k,\, l}$ is the base of the tautological $M_k(\mathbb{C})$-bundle
(its fiber over a point $\alpha \in \Gr_{k,\, l}$ is the $k$-subalgebra $M_{k,\, \alpha}\subset M_{kl}(\mathbb{C})$ corresponding to this point)
which we denote by ${\mathcal A}_{k,\, l}\rightarrow \Gr_{k,\, l}$.
More precisely, ${\mathcal A}_{k,\, l}$ is a subbundle in $\Gr_{k,\, l}\times M_{kl}(\mathbb{C})$ consisting of all pairs
$(\alpha,\, T),\: \alpha \in \Gr_{k,\, l},\, T\in M_{k,\, \alpha},$ where $M_{k,\, \alpha}$ is the $k$-subalgebra in $M_{kl}(\mathbb{C})$
corresponding to the point $\alpha \in \Gr_{k,\, l}$.

Let ${\mathcal B}_{k,\,l}\rightarrow \Gr_{k,\, l}$ be the bundle of centralizers for the subbundle
${\mathcal A}_{k,\, l}\subset \Gr_{k,\, l}\times M_{kl}(\mathbb{C})$. It is easy to see that
${\mathcal B}_{k,\,l}\rightarrow \Gr_{k,\, l}$ has fiber $M_l(\mathbb{C})$.

It follows from representations (\ref{reprfr}) and (\ref{reprmgr}) that $\Gr_{k,\, l}$ is the base of the principal
$\PU(k)$-bundle $\pi_{k,\, l}\colon \Fr_{k,\, l}\rightarrow \Gr_{k,\, l}.$
Clearly, ${\mathcal A}_{k,\, l}\rightarrow \Gr_{k,\, l}$ is associated with this principal bundle
with respect to the action $\PU(k)\stackrel{\cong}{\rightarrow}\Aut(M_k(\mathbb{C}))$
(recall that we consider $*$-automorphisms only).
Hence the pullback $\pi^*_{k,\, l}({\mathcal A}_{k,\, l})$ has the canonical trivialization
(while the bundle $\pi^*_{k,\, l}({\mathcal B}_{k,\, l})\rightarrow \Fr_{k,\, l}$ is nontrivial, see below).

In general, $\mu$ (see (\ref{muu})) is a nontrivial embedding, in particular, it can be nonhomotopic to the choice of a constant
$k$-subalgebra in $X\times M_{kl}(\mathbb{C})$ (in this case the homotopy class of $X\rightarrow \Fr_{k,\, l}$ is nontrivial).
In particular, the subbundle $B_l\rightarrow X$ (with fiber
$M_l(\mathbb{C})$) in $X\times M_{kl}(\mathbb{C})$ of centralizers for $\mu(X\times M_k(\mathbb{C}))\subset
X\times M_{kl}(\mathbb{C})$ can be nontrivial.

The fibration
\begin{equation}
\label{rassl}
\PU(l)\stackrel{E_k\otimes
\ldots}{\longrightarrow}\PU(kl)\stackrel{\chi_k^\prime}{\longrightarrow}\Fr_{k,\, l}
\end{equation}
(cf. (\ref{reprfr})) can be extended to the right
\begin{equation}
\label{fibbrr}
\Fr_{k,\,
l}\stackrel{\psi_k^\prime}{\longrightarrow}\BPU(l)\stackrel{\omega_k^\prime}{\longrightarrow}\BPU(kl),
\end{equation}
where $\psi_k^\prime$ is the classifying map for the
$M_l(\mathbb{C})$-bundle $\widetilde{\mathcal
B}_{k,\,l}:=\pi^*_{k,\, l}({\mathcal B}_{k,\,l})\rightarrow \Fr_{k,\, l}$ (which is associated with the principal $\PU(l)$-bundle (\ref{rassl})).

Let $[M_k]$ be the trivial $M_k(\mathbb{C})$-bundle $X\times M_k(\mathbb{C})$ over $X$.

\begin{proposition}
\label{ktors2}
(Cf. Proposition \ref{ktors})\, For an
$M_l(\mathbb{C})$-bundle $B_l\rightarrow X$ such that
\begin{equation}
\label{condtriv2}
[M_k]\otimes B_l\cong X\times M_{kl}(\mathbb{C})
\end{equation}
(cf. (\ref{condtriv})) a classifying map
$\varphi_{B_l}\colon X\rightarrow \BPU(l)$ can be lifted to
$\widetilde{\varphi}_{B_l}\colon X\rightarrow
\Fr_{k,\, l}$ (i.e. $\psi_k^\prime \circ
\widetilde{\varphi}_{B_l}=\varphi_{B_l}$ or
$B_l=\widetilde{\varphi}_{B_l}^*(\widetilde{\mathcal
B}_{k,\,l})$).
\end{proposition}
{\noindent \it Proof \;} follows from the analysis of fibration (\ref{fibbrr}).$\quad \square$

\smallskip

Moreover, the choice of such a lift corresponds to the choice of trivialization (\ref{condtriv2}) and we
return to the interpretation of the map $X\rightarrow \Fr_{k,\, l}$ given in Proposition \ref{interpret}.
We stress that a map $X\rightarrow \Fr_{k,\, l}$ is not just an
$M_l(\mathbb{C})$-bundle, but an $M_l(\mathbb{C})$-bundle together with a particular choice of
trivialization (\ref{condtriv2}).

It is not difficult to show \cite{Ers2} that the bundle
$B_l\rightarrow X$ as in the statement of Proposition \ref{ktors2} has the form $\End(\eta_l)$
for some (unique up to isomorphism) $\mathbb{C}^l$-bundle $\eta_l\rightarrow X$ with the structural group
$\SU(l)$ (here the condition $(k,\, l)=1$ is essential).

Let $\widetilde{\zeta}\rightarrow
\Fr_{k,\, l}$ be the line bundle associated with the universal covering
$\rho_k\rightarrow \widetilde{\Fr}_{k,\, l}\rightarrow \Fr_{k,\, l}$, where $\rho_k$ is the group of $k$th roots of unity.
Note that $\widetilde{\Fr}_{k,\, l}=\SU(kl)/(E_k\otimes \SU(l))$.
Put $\zeta^\prime :=\widetilde{\varphi}_{B_l}^*(\widetilde{\zeta})\rightarrow X$ and $\eta_l^\prime :=\eta_l\otimes
\zeta^\prime.$

Recall that $\Fred_n({\mathcal H})$ is the subspace of Fredholm operators in
$M_n({\mathcal B}({\mathcal H}))$. The evaluation map
\begin{equation}
\label{canmappp}
ev_{k,\, l}\colon \Fr_{k,\, l}\times M_k(\mathbb{C})\rightarrow
M_{kl}(\mathbb{C}),\quad ev_{k,\, l}(h,\, T)=h(T)
\end{equation}
(recall that $\Fr_{k,\, l}:=\Hom_{alg}(M_k(\mathbb{C}),\,
M_{kl}(\mathbb{C}))$) induces the map
\begin{equation}
\label{actreqloc}
\gamma_{k,\, l}^\prime \colon \Fr_{k,\, l}\times
\Fred_k({\mathcal H})\rightarrow \Fred_{kl}({\mathcal H}).
\end{equation}

\begin{remark}
Note that map (\ref{canmappp}) can be decomposed into the composition
$$
\Fr_{k,\, l}\times M_k(\mathbb{C})\rightarrow \Fr_{k,\,
l}{\mathop{\times}\limits_{\PU(k)}}M_k(\mathbb{C})={\mathcal
A}_{k,\, l}\rightarrow M_{kl}(\mathbb{C}),
$$
where the last map is the tautological embedding $\mu \colon
{\mathcal A}_{k,\, l}\rightarrow \Gr_{k,\, l}\times
M_{kl}(\mathbb{C})$ followed by the projection onto the second factor.
\end{remark}

Let $f_\xi \colon X\rightarrow \Fred_k({\mathcal H})$
represent some element $\xi \in K(X).$

\begin{theorem}
\label{mth}
(Cf. Theorem \ref{lineth}). With respect to the above notation the composite map (cf. (\ref{compmapconj}))
$$
X\stackrel{\diag}{\longrightarrow}X\times
X\stackrel{\widetilde{\varphi}_{B_l} \times
f_\xi}{\longrightarrow}\Fr_{k,\, l}\times \Fred_k({\mathcal H})
\stackrel{\gamma_{k,\,
l}^\prime}{\longrightarrow}\Fred_{kl}({\mathcal H})
$$
represents the element $\eta_l^\prime \otimes \xi \in K(X).$
\end{theorem}
{\noindent \it Proof}\; (cf. \cite{AS}, Proposition 2.1).
By assumption the element $\xi \in K(X)$ is represented by a family of Fredholm operators $F=\{ F_x\}$ in
a Hilbert space ${\mathcal H}^k$. Then the element $\eta_l^\prime \otimes \xi \in K(X)$ is represented by the family
of Fredholm operators $\{ \Id_{(B_l)_x}\otimes \, F_x\}$ in the Hilbert bundle $\eta_l^\prime \otimes ({\mathcal H}^k)$
(recall that $\End(\eta_l)=B_l\, \Rightarrow \, \End(\eta_l^\prime)=B_l$). A trivialization
$\eta_l^\prime \otimes ({\mathcal H}^k)\cong {\mathcal H}^{kl}$ is the same thing as a map $\widetilde{\varphi}_{B_l}\colon X\rightarrow \Fr_{k,\, l}$, i.e. a lift
of the classifying map $\varphi_{B_l}\colon X\rightarrow \BPU(l)$ for $B_l$ (see (\ref{fibbrr})).$\quad \square$

\smallskip

\begin{remark}
In order to separate the ``$\SU$''-part of the ``action'' $\gamma_{k,\, l}^\prime$ from its ``line'' part,
one can use the space $\widetilde{\Fr}_{k,\, l}=\SU(kl)/(E_k\otimes \SU(l))$ \cite{Ers2} in place of $\Fr_{k,\, l}$.
Then one would have the representing map for $\eta_l \otimes \xi \in K(X)$ instead of $\eta_l^\prime \otimes \xi$ in the statement of Theorem \ref{mth}.
\end{remark}

\begin{remark}
Note that $\Fr_{k,\, 1}=\PU(k)$ and the action $\gamma_{k,\, 1}^\prime$ coincides with the action $\gamma_{k}^\prime$
from the previous section.
\end{remark}

Now using the composition of algebra homomorphisms we are going to define maps
$\phi_{k,\, l} \colon \Fr_{k,\, l}\times \Fr_{k,\, l}\rightarrow \Fr_{k,\, l^2},$ i.e.
$$
\phi_{k,\, l} \colon \Hom_{alg}(M_k(\mathbb{C}),\, M_{kl}(\mathbb{C}))\times \Hom_{alg}(M_k(\mathbb{C}),\, M_{kl}(\mathbb{C}))
\rightarrow \Hom_{alg}(M_k(\mathbb{C}),\, M_{kl^2}(\mathbb{C})).
$$
First let us define a map
$$
\iota_{k,\, l} \colon \Hom_{alg}(M_k(\mathbb{C}),\, M_{kl}(\mathbb{C}))\rightarrow \Hom_{alg}(M_{kl}(\mathbb{C}),\, M_{kl^2}(\mathbb{C})),\;
\iota_{k,\, l}(h)=h\otimes \id_{M_l(\mathbb{C})}.
$$
Then $\phi_{k,\, l}$ is defined as the composition of homomorphisms: $\phi_{k,\, l}(h_2,\, h_1)=\iota_{k,\, l}(h_2)\circ h_1,$
where $h_i\in \Hom_{alg}(M_k(\mathbb{C}),\, M_{kl}(\mathbb{C})).$
Then we have $ev_{k,\, l^2}(\phi_{k,\, l}(h_2,\, h_1),\, T)=ev_{kl,\, l}(\iota_{k,\, l}(h_2),\, ev_{k,\, l}(h_1,\, T)),$
i.e. $\phi_{k,\, l}(h_2,\, h_1)(T)=\iota_{k,\, l}(h_2)(h_1(T)),$ where $T\in M_k(\mathbb{C})$.

Now suppose there is an $M_l(\mathbb{C})$-bundle $C_l\rightarrow X$ with the corresponding vector bundle $\theta_l^\prime$
such that $C_l\cong \End(\theta_k^\prime)$ (cf. a few paragraphs after Proposition \ref{ktors2}). Suppose that $\widetilde{\varphi}_{C_l}\colon X\rightarrow \Fr_{k,\, l}$ is
its classifying map.

\begin{proposition}
(cf. Theorem \ref{mth}.) The composition
$$
X\stackrel{\diag}{\rightarrow}X\times X\times X\stackrel{\widetilde{\varphi}_{C_l}\times \widetilde{\varphi}_{B_l}\times f_\xi}{\longrightarrow}
\Fr_{k,\, l}\times \Fr_{k,\, l}\times \Fred_k({\mathcal H})\rightarrow \Fred_{kl^2}({\mathcal H})
$$
where the last map is the composition $\gamma^\prime_{k,\, l^2}\circ (\phi_{k,\, l}\times \id_{\Fred_k({\mathcal H})})=
\gamma^\prime_{kl,\, l}\circ (\iota_{k,\, l}\times \gamma^\prime_{k,\, l})$
represents the element $\eta_l^\prime \otimes \theta_l^\prime\otimes \xi \in K(X).$
\end{proposition}
{\noindent \it Proof}\; is evident.$\quad \square$

\smallskip

Clearly, the results of this section can be generalized to the case of spaces $\Fr_{kl^m,\, l^n},\; m,\, n\in \mathbb{N}.$
In the next section we will construct a genuine action of their direct limit $\Fr_{kl^\infty,\, l^\infty}$
on the space of Fredholm operators.

\section{A construction of the classifying space}

A simple calculation with homotopy groups shows that the direct limit
$\Fr_{kl^\infty,\, l^\infty}:=\varinjlim_{n}\Fr_{kl^n,\,
l^n}$ for $(k,\, l)=1$ is not contractible because its homotopy groups are $\mathbb{Z}/k\mathbb{Z}$ in odd dimensions (and $0$ in even ones).
In this section we
show that it is a topological monoid and construct its
classifying space.

\subsection{The category ${\mathcal C}_{k,\, l}$}

First, we define some auxiliary category ${\mathcal C}_{k,\, l}$.
Fix a pair of positive integers $\{ k,\, l\} ,\; k,\, l>1, (k,\, l)=1.$
By ${\mathcal C}_{k,\, l}$ denote the category with the countable number of objects which are
matrix algebras of the form
$M_{kl^m}(\mathbb{C})$ ($m=0,\, 1,\, \ldots)$ and morphisms in ${\mathcal C}_{k,\, l}$ from
$M_{kl^m}(\mathbb{C})$ to $M_{kl^n}(\mathbb{C})$ are all unital
$*$-homomorphisms $M_{kl^m}(\mathbb{C})\rightarrow M_{kl^n}(\mathbb{C})$ (this set is nonempty iff $m\leq n$).
(Since $k>1$ we see that this category does not contain the initial object, therefore there is no reason to expect that
its classifying space (i.e. the ``geometric realization'') is contractible).
Thus, morphisms
$\Mor(M_{kl^m}(\mathbb{C}),\: M_{kl^n}(\mathbb{C}))$
form the space $\Fr_{kl^m,\, l^{n-m}}=\Hom_{alg}(M_{kl^m}(\mathbb{C}),\,
M_{kl^n}(\mathbb{C}))$, hence ${\mathcal C}_{k,\, l}$ is a topological category.
(Note that for $n=m$ the space $\Fr_{kl^m,\, l^{n-m}}=\Fr_{kl^m,\, 1}$ is the group $\PU(kl^m).$)
In particular, there is the collection of continuous maps
$\Fr_{kl^{m+n},\, l^{r}}\times \Fr_{kl^m,\, l^{n}}\rightarrow \Fr_{kl^m,\, l^{n+r}}$
for all $m,\, n,\, r\geq 0$ given by the composition of morphisms.

Recall (\cite{Weibel}, Chapter IV) that there is an appropriate modification of the construction of the
geometric realization
$\B{\mathcal C}$ for topological categories ${\mathcal C}$.
More precisely, the nerve in this case is a simplicial {\it topological space}
and $\B{\mathcal C}$ is its appropriate geometric realization.
Now we are going to describe the classifying space of the category $\B {\mathcal C}_{k,\, l}$.

So let $\B {\mathcal C}_{k,\, l}$ be the classifying space of the topological category
${\mathcal C}_{k,\, l}$. Its $0$-cells (vertices) are objects of ${\mathcal C}_{k,\, l}$,
i.e. actually positive integers.
Its $1$-cells (edges) are morphisms in
${\mathcal C}_{k,\, l}$ (excluding identity morphisms) attached to their
source and target, i.e. $\coprod_{m,\, n\geq 0} \Fr_{kl^m,\, l^n}$ (recall that for $n=0\; \Fr_{kl^m,\, 1}=\PU(kl^m)$).
For each pair of composable morphisms
$h_0,\, h_1$ in ${\mathcal C}_{k,\, l}$ there is a $2$-simplex:
$$
\diagram
& 1 \drto^{h_1} \\
0 \rrto_{h_1\circ h_0} \urto^{h_0} && 2 \\
\enddiagram
$$
attached to the $1$-skeleton, etc. Thus, at the third step we obtain the space
$\coprod_{m,\, n,\, r\geq 0} (\Fr_{kl^{m+n},\, l^r}\times \Fr_{kl^m,\, l^n})$
consisting of all pairs of composable morphisms.
The face maps $\partial_0,\, \partial_2$ are defined by the deletion of the corresponding morphism ($h_0$ and $h_1$ respectively in the above simplex) and
$\partial_1$ is defined by the composition map $\coprod_{m,\, n,\, r\geq 0} (\Fr_{kl^{m+n},\, l^r}\times \Fr_{kl^m,\, l^n})\rightarrow
\coprod_{m,\, n+r\geq 0} \Fr_{kl^m,\, l^{n+r}}=\coprod_{m,\, n\geq 0} \Fr_{kl^m,\, l^n}.$
The nerve $N{\mathcal C}_{k,\, l}$ of the category ${\mathcal C}_{k,\, l}$ is
\begin{equation}
\label{nervC}
(\mathbb{N},\, \coprod_{m,\, n\geq 0} \Fr_{kl^m,\, l^n},\, \coprod_{m,\, n,\, r\geq 0} (\Fr_{kl^{m+n},\, l^r}\times \Fr_{kl^m,\, l^n}),\ldots ).
\end{equation}

Recall that there is a construction of the classifying space of a (topological) group $G$
as the geometric realization of the simplicial topological space
$(\pt,\, G,\, G\times G,\ldots )$. Then the total space $\EG$ of
the universal principal $G$-bundle is the geometric realization of the simplicial space $(G,\, G\times G,\ldots )$.
Consider the simplicial topological space
$$
{\mathcal E}_{k,\, l}:=
(\coprod_{m,\, n\geq 0} \Fr_{kl^m,\, l^n},\, \coprod_{m,\, n,\, r\geq 0} (\Fr_{kl^{m+n},\, l^r}\times \Fr_{kl^m,\, l^n}),\ldots ),
$$
whose faces and degeneracies are defined by analogy with the construction of $\EG$.
There is the map of simplicial spaces $p\colon {\mathcal E}_{k,\, l}\rightarrow N{\mathcal C}_{k,\, l}.$

Now we are going to give a construction of the corresponding ``universal bundle''.
The idea is to construct a ``simplicial bundle'' associated with the ``universal principal bundle'' $p$.

More precisely, again consider the space ${\mathcal E}_{k,\, l}$:
$$
(\coprod_{m,\, n\geq 0} \Fr_{kl^m,\, l^n},\, \coprod_{m,\, n,\, r\geq 0} (\Fr_{kl^{m+n},\, l^r}\times \Fr_{kl^m,\, l^n}),\,
\coprod_{m,\, n,\, r,\, s\geq 0} (\Fr_{kl^{m+n+r},\, l^s}\times \Fr_{kl^{m+n},\, l^r}\times \Fr_{kl^m,\, l^n}),\ldots ).
$$
Applying the natural (``evaluation'') maps $\Fr_{kl^m,\, l^n}\times M_{kl^m}(\mathbb{C})\rightarrow M_{kl^{m+n}}(\mathbb{C})$, we obtain
\begin{equation}
\label{nervB}
(\coprod_{m\geq 0} M_{kl^{m}}(\mathbb{C}),\, \coprod_{m,\, n\geq 0} (\Fr_{kl^{m},\, l^n}\times M_{kl^{m}}(\mathbb{C})),\,
\coprod_{m,\, n,\, r\geq 0} (\Fr_{kl^{m+n},\, l^r}\times \Fr_{kl^{m},\, l^n}\times M_{kl^{m}}(\mathbb{C})),\ldots ).
\end{equation}
The obtained object can be regarded as a ``simplicial bundle'' over (\ref{nervC}).
Now performing the appropriate factorizations we should define its geometric realization. Namely, the matrix algebras from the first disjoint
union in (\ref{nervB}) are fibers of our bundle over $0$-cells of the space $\B {\mathcal C}_{k,\, l}$, from the second over $1$-cells, etc.
The attachment of $1$-cells to their source corresponds to the projection $\Fr_{kl^{m},\, l^n}\times M_{kl^{m}}(\mathbb{C})\rightarrow M_{kl^{m}}(\mathbb{C})$
onto the second factor, and the one to their target corresponds to the natural map $\Fr_{kl^m,\, l^n}\times M_{kl^m}(\mathbb{C})\rightarrow M_{kl^{m+n}}(\mathbb{C})$.
$M_{kl^{m}}(\mathbb{C})$-bundles over $2$-simplices correspond to the products
$\Fr_{kl^{m+n},\, l^r}\times \Fr_{kl^{m},\, l^n}\times M_{kl^{m}}(\mathbb{C})$ and they are identified over the boundary as follows:
$\overline{\partial}_0\colon \Fr_{kl^{m+n},\, l^r}\times \Fr_{kl^{m},\, l^n}\times M_{kl^{m}}(\mathbb{C})\rightarrow
\Fr_{kl^{m+n},\, l^r}\times M_{kl^{m+n}}(\mathbb{C})$ is induced by the natural map $\Fr_{kl^m,\, l^n}\times M_{kl^m}(\mathbb{C})\rightarrow M_{kl^{m+n}}(\mathbb{C})$,
$\overline{\partial}_1\colon \Fr_{kl^{m+n},\, l^r}\times \Fr_{kl^{m},\, l^n}\times M_{kl^{m}}(\mathbb{C})\rightarrow
\Fr_{kl^{m},\, l^{n+r}}\times M_{kl^{m}}(\mathbb{C})$ by the composition of morphisms $\Fr_{kl^{m+n},\, l^r}\times \Fr_{kl^{m},\, l^n}\rightarrow
\Fr_{kl^{m},\, l^{n+r}}$, а $\overline{\partial}_2\colon \Fr_{kl^{m+n},\, l^r}\times \Fr_{kl^{m},\, l^n}\times M_{kl^{m}}(\mathbb{C})\rightarrow
\Fr_{kl^{m},\, l^n}\times M_{kl^{m}}(\mathbb{C})$ by the projection onto the product of the second and the third factors.
For bundles over simplices of greater dimension the attaching maps are defined analogously.

The obtained topological space has the natural projection onto
the space $\B {\mathcal C}_{k,\, l}$ with fibres $M_{kl^{m}}(\mathbb{C})$.

Note that in the same way as ${\mathcal C}_{k,\, l}$ one can define categories ${\mathcal C}_{k^m,\, l},\; m\geq 1.$
The objects of ${\mathcal C}_{k^m,\, l}$ are matrix algebras $M_{k^ml^n}(\mathbb{C}),\; n\geq 0$.
Note that the tensor product of matrix algebras gives rise to the bifunctor
$T_{m,\, n}\colon {\mathcal C}_{k^m,\, l}\times {\mathcal C}_{k^n,\, l}\rightarrow {\mathcal C}_{k^{m+n},\, l}$.
Thus, on objects we have:
$T_{m,\, n}(M_{k^ml^r}(\mathbb{C}),\, M_{k^nl^s}(\mathbb{C}))=M_{k^ml^r}(\mathbb{C})\otimes M_{k^nl^s}(\mathbb{C})=M_{k^{m+n}l^{r+s}}(\mathbb{C})$,
and on morphisms $h_1\colon M_{k^ml^r}(\mathbb{C})\rightarrow M_{k^ml^{r+t}}(\mathbb{C}),\; h_2\colon M_{k^nl^s}(\mathbb{C})\rightarrow M_{k^nl^{s+u}}(\mathbb{C})\quad
T_{m,\, n}(h_1,\, h_2)=h_1\otimes h_2\colon M_{k^{m+n}l^{r+s}}(\mathbb{C})\rightarrow M_{k^{m+n}l^{r+s+t+u}}(\mathbb{C}).$
Moreover, the bifunctor $T_{m,\, n}$ determines the continuous map
of the topological spaces $\Fr_{k^ml^r,\, l^t}\times \Fr_{k^nl^s,\, l^u}\rightarrow \Fr_{k^{m+n},\, l^{r+s+t+u}},\; (h_1,\, h_2)\mapsto h_1\otimes h_2$.

For the topological bicategory
${\mathcal C}_{k^m,\, l}\times {\mathcal C}_{k^n,\, l}$ one can define the bisimplicial topological space
$X=\{ X_{p,\, q}\} ,$ where $X_{p,\, q}$ consists of all pairs of functors ${\bf p+1}\rightarrow {\mathcal C}_{k^m,\, l},\;
{\bf q+1}\rightarrow {\mathcal C}_{k^n,\, l}.$ There are two types of face and degeneracy maps: ``horizontal'' and ``vertical''
which correspond to the category ${\mathcal C}_{k^m,\, l}$ and
${\mathcal C}_{k^n,\, l}$ respectively. Clearly that its geometric realization $\B X$ (to any point of $X_{p,\, q}$
we attach $\Delta^p\times \Delta^q$) is
$\B {\mathcal C}_{k^m,\, l}\times \B {\mathcal C}_{k^n,\, l}.$
The bifunctor $T_{m,\, n}$ determines the continuous map
$\widetilde{T}_{m,\, n}\colon \B {\mathcal C}_{k^m,\, l}\times \B {\mathcal C}_{k^n,\, l}\rightarrow \B {\mathcal C}_{k^{m+n},\, l}.$
It seems that the family of such bifunctors with different
$m,\, n$ defines a structure of $H$-space on the direct limit $\varinjlim_m\B {\mathcal C}_{k^m,\, l}$
(the maps $\B {\mathcal C}_{k^m,\, l}\rightarrow \B {\mathcal C}_{k^n,\, l}$ are given by functors
${\mathcal C}_{k^m,\, l}\rightarrow {\mathcal C}_{k^n,\, l}$ which on objects are defined by the tensor product of matrix algebras
$M_{k^ml^r}(\mathbb{C}),\; r\geq 0$
by the fixed $M_{k^{n-m}}(\mathbb{C})$ and on morphisms by the assignment $h\mapsto h\otimes \id_{M_{k^{n-m}}(\mathbb{C})}$).

In conclusion of this subsection we describe yet another two properties of the category ${\mathcal C}_{k,\, l}$.

Let ${\mathcal C}_{k,\, l}^a$
be the category with the same objects as
${\mathcal C}_{k,\, l}$ but with morphisms that are automorphisms in
${\mathcal C}_{k,\, l}$. Clearly,
$\B {\mathcal C}_{k,\, l}^a\simeq \coprod_{n\geq 0}\BPU(kl^n)$ and
the embedding $\coprod_{n\geq 0}\BPU(kl^n)\rightarrow \B {\mathcal C}_{k,\, l}$ corresponds to
the inclusion of the subcategory ${\mathcal C}_{k,\, l}^a\rightarrow {\mathcal C}_{k,\, l}$.

Let ${\bf N}$ be the category with countable set of objects $\{ 0,\, 1,\, 2,\, \ldots \}$ and there is a morphism
from $i$ to $j$ iff $i\leq j$, and such morphism is unique.
It is easy to see that its classifying space $\B {\bf N}$
is the infinite simplex $\Delta$.

There is the obvious functor $F\colon {\mathcal C}_{k,\, l}\rightarrow {\bf N},\; F(M_{kl^m}(\mathbb{C}))=m$ and for a morphism
$h\colon M_{kl^m}(\mathbb{C})\rightarrow M_{kl^n}(\mathbb{C})$ the morphism $F(h)$
is the unique morphism $m\rightarrow n$ in ${\bf N}$. Therefore there is the corresponding map of classifying spaces
$\B {\mathcal C}_{k,\, l}\rightarrow \B {\bf N}.$
The subspace $\coprod_{n\geq 0}\BPU(kl^n)\subset \B {\mathcal C}_{k,\, l}$
corresponds to the vertices of the simplex $\B {\bf N}$ (more precisely, to the corresponding discrete category).
In general, simplices degenerate under this map of classifying spaces,
moreover, their degenerations correspond to automorphisms of objects of
${\mathcal C}_{k,\, l}$ and nondegenerate simplices correspond to chains
$\{ h_0,\, h_1,\, \ldots h_n\}$ of composable morphisms such that for all $h_r$ the source and the target are different objects
(in other words, $h_i\in \Fr_{kl^m,\, l^n},\; n\neq 0$).

\subsection{The definition of $\B \Fr_{kl^\infty,\, l^\infty}$}

Now consider new topological category $\overline{\mathcal C}_{k,\, l}$. It has a unique object $M_{kl^\infty}(\mathbb{C})=\varinjlim_mM_{kl^m}(\mathbb{C})$
(i.e. $\overline{\mathcal C}_{k,\, l}$ is actually a monoid), where the direct limit is taken over unital $*$-homomorphisms $M_{kl^m}(\mathbb{C})\rightarrow M_{kl^{m+1}}(\mathbb{C}),\;
X\mapsto X\otimes E_l.$ More precisely, we assume that the matrix algebra $M_{kl^\infty}(\mathbb{C})$ is given together with the infinite family
of $*$-subalgebras
$A_k\subset A_{kl}\subset A_{kl^2}\subset \ldots,$ where $A_{kl^m}=M_{kl^m}(\mathbb{C})$ and $A_{kl^{m+1}}=M_l(A_{kl^m})$ for every $m\geq 0$,
which form its filtration.

By definition, for each morphism $h\colon M_{kl^\infty}(\mathbb{C})\rightarrow M_{kl^\infty}(\mathbb{C}),\; h\in \Mor(\overline{\mathcal C}_{k,\, l})$
there exists a pair $m,\, n,\: n\geq m$ such that 1) $h|_{A_{kl^m}}$ is a unital $*$-homomorphism $h|_{A_{kl^m}}\colon A_{kl^m}\rightarrow A_{kl^n}$
and 2) $h=M_{l^\infty}(h|_{A_{kl^m}})$, i.e.
$h|_{A_{kl^{m+1}}}=M_l(h|_{A_{kl^m}})\colon A_{kl^{m+1}}=M_l(A_{kl^m})\rightarrow A_{kl^{n+1}}=M_l(A_{kl^n}),$ etc.
In other words, $h$ is induced by some $h^\prime \in \Fr_{kl^m,\, l^{n-m}}=\Hom_{alg}(A_{kl^m},\, A_{kl^n})$.
For example, the composition of morphisms induced by $h_1\colon A_{kl}\rightarrow A_{kl^2}$
and $h_2\colon A_k\rightarrow A_{kl}$ is defined by the following diagram
$$
\diagram
\ldots & \ldots & \ldots \\
A_{kl^4} \uto^{\cup} & A_{kl^4} \uto^{\cup} & A_{kl^4} \uto^{\cup} \\
A_{kl^3} \uto^{\cup} \urto^{\qquad \qquad \qquad \qquad M_{l^2}(h)} & A_{kl^3} \uto^{\cup} \urto^{\qquad \qquad \qquad \qquad M_{l^3}(h_2)} & A_{kl^3} \uto^{\cup} \\
A_{kl^2} \uto^{\cup} \urto^{\qquad \qquad \qquad \qquad M_l(h_1)} & A_{kl^2} \uto^{\cup} \urto^{\qquad \qquad \qquad \qquad M_{l^2}(h_2)} & A_{kl^2} \uto^{\cup} \\
A_{kl} \uto^{\cup} \urto^{h_1} & A_{kl} \uto^{\cup} \urto^{\qquad \qquad \qquad \qquad M_l(h_2)} & A_{kl} \uto^{\cup} \\
A_{k} \uto^{\cup} & A_{k} \uto^{\cup} \urto^{h_2} & A_{k} \uto^{\cup} \\
\enddiagram
$$
as the class of the homomorphism $M_{l^2}(h_2)\circ h_1.$
Clearly, the composition of morphisms is well-defined and associative and the identity morphism is $M_{l^\infty}(\id_{A_k}),$
i.e. the family $\{ \id_{A_k},\, \id_{A_{kl}},\, \id_{A_{kl^2}},\, \ldots \,\}$.

Now we define the functor $\Phi \colon {\mathcal C}_{k,\, l}\rightarrow \overline{\mathcal C}_{k,\, l}$
which sends every object $M_{kl^m}(\mathbb{C})$ in ${\mathcal C}_{k,\, l}$ to the unique object $M_{kl^\infty}(\mathbb{C})$ in $\overline{\mathcal C}_{k,\, l}$.
For a morphism $h\in \Fr_{kl^m,\, l^{n-m}}$ we put $\Phi(h)=M_{l^\infty}(h).$

Thus, we see that $\Mor(\overline{\mathcal C}_{k,\, l})$ is the {\bf well-pointed grouplike}
(because $\pi_0(\Fr_{kl^\infty,\, l^\infty})=0$) {\bf topological monoid} $\Fr_{kl^\infty,\, l^\infty}.$
Recall \cite{Rudyak} that for such a monoid $M$ there exists the classifying space $\B M$. Thus we have
{\bf the classifying space} $\B \Fr_{kl^\infty,\, l^\infty}$ which is defined uniquely up to $CW$-equivalence
and there is the Whitehead equivalence $\Fr_{kl^\infty,\, l^\infty}\rightarrow \Omega \B \Fr_{kl^\infty,\, l^\infty}$;
in particular, $\pi_i(\Fr_{kl^\infty,\, l^\infty})=\pi_{i+1}(\B \Fr_{kl^\infty,\, l^\infty})$.
Note that $\Phi$ defines a continuous map $\widetilde{\Phi}\colon \B {\mathcal C}_{k,\, l}\rightarrow \B \Fr_{kl^\infty,\, l^\infty}.$
Moreover, the maps $\widetilde{T}_{m,\, n}$ correspond to the maps $\B \Fr_{k^ml^\infty,\, l^\infty}\times \B \Fr_{k^nl^\infty,\, l^\infty}
\rightarrow \B \Fr_{k^{m+n}l^\infty,\, l^\infty}$ which are given by maps $\Fr_{k^ml^\infty,\, l^\infty}\times \Fr_{k^nl^\infty,\, l^\infty}
\rightarrow \Fr_{k^{m+n}l^\infty,\, l^\infty}$ induced by the tensor product of matrix algebras.

Note that there is the subgroup $\PU(kl^\infty)=\varinjlim_m\PU(kl^m)$ in the monoid $\Fr_{kl^\infty,\, l^\infty}$ which corresponds to automorphisms
of $M_{kl^\infty}(\mathbb{C})$.
The corresponding map $\BPU(kl^\infty)\rightarrow \B \Fr_{kl^\infty,\, l^\infty}$ has the homotopy fiber $\Gr_{kl^\infty,\, l^\infty}$
(this follows from the fibration
$$
\PU(kl^m)\rightarrow \Fr_{kl^m,\, l^n}=\PU(kl^{m+n})/(E_{kl^m}\otimes \PU(l^n))\rightarrow \Gr_{kl^m,\, l^n}).
$$

\begin{remark}
Retrospectively, we note that we have used the maps $\Fr_{kl^{m+n},\, l^r}\times \Fr_{kl^m,\, l^n}\rightarrow \Fr_{kl^m,\, l^{n+r}}$
(given by the composition of morphisms) in order to define the monoid structure on $\Fr_{kl^\infty,\, l^\infty}$
and thereby the classifying space $\B \Fr_{kl^\infty,\, l^\infty}$, and the maps $\Fr_{k^ml^r,\, l^s}\times \Fr_{k^nl^t,\, l^u}\rightarrow
\Fr_{k^{m+n}l^{r+t},\, l^{l^{s+u}}}$ given by the tensor product in order to define the additional structure
on the spaces $\Fr_{k^ml^\infty,\, l^\infty}$ (which gives rise to the $H$-space structure on $\varinjlim_m\Fr_{k^ml^\infty,\, l^\infty}$).
From the category-theoretic point of view the first corresponds to the composition of morphisms in the category and the second
to the monoidal structure on it.
\end{remark}

\subsection{The action of $\Fr_{kl^\infty,\, l^\infty}$ on the space of Fredholm operators}

Recall (\ref{canmappp}) that there are evaluation maps
$$
ev_{kl^m,\, l^{n-m}}\colon \Fr_{kl^m,\, l^{n-m}}\times M_{kl^m}(\mathbb{C})\rightarrow M_{kl^n}(\mathbb{C})
$$
and the corresponding maps (recall that $\Fred_{kl^m}({\mathcal H})$ is the subspace of Fredholm operators in $M_{kl^m}({\mathcal B}({\mathcal H}))$)
$$
\gamma^\prime_{kl^m,\, l^{n-m}}\colon \Fr_{kl^m,\, l^{n-m}}\times \Fred_{kl^m}({\mathcal H})\rightarrow \Fred_{kl^n}({\mathcal H})
$$
(see (\ref{actreqloc})).
Using the filtration in $M_{kl^\infty}({\mathcal B}({\mathcal H}))$ (and hence in $\varinjlim_m\Fred_{kl^m}({\mathcal H})$) corresponding to the above filtration
$A_k\subset A_{kl}\subset A_{kl^2}\subset \ldots$ in the matrix algebra $M_{kl^\infty}(\mathbb{C})$ one can define the action of the monoid
$\Fr_{kl^\infty,\, l^\infty}$ on $\varinjlim_m\Fred_{kl^m}({\mathcal H})$. Note that since the direct limit is taken over maps induced by the tensor product,
we see that $\Fred_{kl^\infty}({\mathcal H}):=\varinjlim_m\Fred_{kl^m}({\mathcal H})$ is the localization in which $l$ becomes invertible
(in particular, the index takes values in $\mathbb{Z}[\frac{1}{l}]$, not in $\mathbb{Z}$).
Thereby we have defined the required action
\begin{equation}
\label{mainaction}
\gamma^\prime_{kl^\infty,\, l^\infty}\colon \Fr_{kl^\infty,\, l^\infty}\times \Fred_{kl^\infty}({\mathcal H})\rightarrow \Fred_{kl^\infty}({\mathcal H})
\end{equation}
of the monoid $\Fr_{kl^\infty,\, l^\infty}$.

Note that the action (\ref{mainaction}) gives rise to the action on $K$-theory (which is recall represented by the space of Fredholm operators)
which corresponds to the action of the $k$-torsion subgroup in $\BU_\otimes$ by tensor products (cf. Proposition \ref{ktors3} and Theorem \ref{mth}).
In fact, this action is defined on $K$-theory $K[\frac{1}{l}]$ localized over $l$ (in the sense that $l$ becomes invertible).
This is not surprising because in (\ref{actreqloc}) we take the tensor product
of $K(X)$ by some $l$-dimensional bundle, $l>1$. It is not difficult to show that in fact our construction does not depend on the choice of $l,\, (k,\, l)=1$.

Note that the restriction of the action $\gamma^\prime_{kl^\infty,\, l^\infty}$ on $\Fr_{kl^m,\, 1}\cong \PU(kl^m)$
coincides with the composition of the action $\gamma^\prime_{kl^m}$ (see Remark \ref{projgract}) and the localization map
$\Fred_{kl^m}({\mathcal H})\rightarrow \Fred_{kl^\infty}({\mathcal H})$ on $l$.

Using this action (\ref{mainaction}) we can define the $\Fred_{kl^\infty}({\mathcal H})$-bundle
$$
\diagram
\Fred_{kl^\infty}({\mathcal H}) \rto & {\rm EFr}_{kl^\infty,\, l^\infty}{\mathop{\times}\limits_{\Fr_{kl^\infty,\, l^\infty}}}\Fred_{kl^\infty}({\mathcal H}) \dto \\
& {\rm BFr}_{kl^\infty,\, l^\infty} \\
\enddiagram
$$
``associated'' with the universal principal $\Fr_{kl^\infty ,\, l^\infty}$-bundle (more precisely,
with the universal principal quasi-fibration, see \cite{Rudyak}) over $\B \Fr_{kl^\infty ,\, l^\infty}$.
This allows us to define a more general version of the twisted $K$-theory than the one given by the action of $Pic(X)$ on $K(X)$.
The above defined maps $\B \Fr_{k^ml^\infty,\, l^\infty}\times \B \Fr_{k^nl^\infty,\, l^\infty}
\rightarrow \B \Fr_{k^{m+n}l^\infty,\, l^\infty}$ give rise to the operation on twistings which is an analog
of the one induced by maps $\BPU(k^m)\times \BPU(k^n)\rightarrow \BPU(k^{m+n})$
in the Abelian case (i.e. the Brauer group), etc.

\section{Appendix: $H$-space $\Fr_{k^\infty,\, l^\infty}$}

In this section we give a category-theoretic description of the structure of $H$-space on $\Fr_{k^\infty,\, l^\infty}$ and $\Gr_{k^\infty,\, l^\infty}$.

\subsection{$k^m$-frames}

\begin{definition}
\label{defframe}
{\it A $k^m$-frame $\alpha$} in the algebra $M_{k^ml^n}(\mathbb{C})$ is an ordered
collection of $k^{2m}$ linearly independent matrices $\{ \alpha_{i,\, j}\}_{1\leq i,\, j\leq k^m}$
such that
\begin{itemize}
\item[(i)] $\alpha_{i,\, j}\alpha_{r,\, s}=\delta_{j,\, r}\alpha_{i,\, s}$ for all $1\leq i,\, j,\, r,\, s\leq k^m$;
\item[(ii)] $\sum_{i=1}^{k^m}\alpha_{i,\, i}=E,$ where $E=E_{k^ml^n}$ is the unit $k^ml^n\times k^ml^n$-matrix which is the unit of the algebra $M_{k^ml^n}(\mathbb{C})$;
\item[(iii)] matrices $\{ \alpha_{i,\, j}\}$ form an orthonormal basis with respect to the hermitian inner product
$(x,\, y):=\tr(x\overline{y}^t)$ on $M_{k^ml^n}(\mathbb{C})$.
\end{itemize}
\end{definition}

For instance, the collection of ``matrix units'' $\{ e_{i,\, j}\}_{1\leq i,\, j\leq k^m}$ (where $e_{i,\, j}$ is the $k^m\times k^m$-matrix whose only nonzero element
is 1 on the intersection of ith row with jth column)
is a $k^m$-frame in $M_{k^m}(\mathbb{C})$,
and the collection $\{ e_{i,\, j}\otimes E_{l^n}\}_{1\leq i,\, j\leq k^m}$ is a $k^m$-frame in $M_{k^ml^n}(\mathbb{C})$.
Clearly, every $k^m$-frame in $M_{k^ml^n}(\mathbb{C})$ is a linear basis in some $k^m$-subalgebra.

\begin{proposition}
The set of all $k^m$-frames in $M_{k^ml^n}(\mathbb{C})$ is the homogeneous space $\PU(k^ml^n)/(E_{k^m}\otimes \PU(l^n)).$
\end{proposition}
{\noindent \it Proof\;} follows from two facts: 1) the group $\PU(k^ml^n)$ of $*$-automorphisms of the algebra $M_{k^ml^n}(\mathbb{C})$
acts transitively on the set of $k^m$-frames,
and 2) the stabilizer of the $k^m$-frame $\{ e_{i,\, j}\otimes E_{l^n}\}_{1\leq i,\, j\leq k^m}$
is the subgroup $E_{k^m}\otimes \PU(l^n)\subset \PU(k^ml^n).\quad \square$

\smallskip

In fact, the space of $k^m$-frames $\Fr_{k^m,\, l^n}$ in $M_{k^ml^n}(\mathbb{C})$ is isomorphic to the space of unital
$*$-homomorphisms
$\Hom_{alg}(M_{k^m}(\mathbb{C}),\, M_{k^ml^n}(\mathbb{C})).$ More precisely, let $\{ e_{i,\, j}\}_{1\leq i,\, j\leq k^m}$ be the frame in $M_{k^m}(\mathbb{C})$ consisting of
matrix units. Then the isomorphism $\Fr_{k^m,\, l^n}\cong \Hom_{alg}(M_{k^m}(\mathbb{C}),\, M_{k^ml^n}(\mathbb{C}))$
is given by the assignment
$$
\alpha \mapsto h_\alpha \colon M_{k^m}(\mathbb{C})\rightarrow M_{k^ml^n}(\mathbb{C}),\; (h_{\alpha})_*(\{ e_{i,\, j}\} )=\alpha \quad \forall \alpha \in \Fr_{k^m,\, l^n}.
$$

Let $\beta$ be a $k^r$-frame in $M_{k^rl^s}(\mathbb{C})$ and $m\leq r.$
Then one can associate with $\beta$ some new $k^m$-frame $\alpha:= \pi_1^m(\beta)$ as follows:
$$
\alpha_{i,\, j}=\beta_{(i-1)k^{r-m}+1,\, (j-1)k^{r-m}+1}+\beta_{(i-1)k^{r-m}+2,\, (j-1)k^{r-m}+2}+\ldots +\beta_{ik^{r-m},\, jk^{r-m}},\quad 1\leq i,\, j \leq k^m.
$$
Also one can associate with $\beta$ some $k^{r-m}$-frame $\gamma:=\pi_2^{r-m}(\beta)$ by the following rule:
$$
\gamma_{i,\, j}=\beta_{i,\, j}+\beta_{i+k^{r-m},\, j+k^{r-m}}+\ldots +\beta_{i+(k^m-1)k^{r-m},\, j+(k^m-1)k^{r-m}},\quad 1\leq i,\, j\leq k^{r-m}.
$$
The idea of the definition of $\pi_1^m(\beta)$ and $\pi_2^{r-m}(\beta)$ is the following. If one takes the $k^r$-frame
$\epsilon$ in $M_{k^r}(\mathbb{C})=M_{k^m}(\mathbb{C})\otimes M_{k^{r-m}}(\mathbb{C})$ consisting of the matrix units, then
the $k^m$ and $k^{r-m}$-frames in subalgebras $M_{k^m}(\mathbb{C})\otimes \mathbb{C}E_{k^{r-m}}\subset M_{k^r}(\mathbb{C})$ and
$\mathbb{C}E_{k^m}\otimes M_{k^{r-m}}(\mathbb{C})\subset M_{k^r}(\mathbb{C})$ consisting of the matrix units tensored by the corresponding
unit matrices are $\pi_1^m(\epsilon)$ and $\pi_2^{r-m}(\epsilon)$ respectively. From the other hand it is easy to see that the frame
$\epsilon$ (under the appropriate ordering) is the tensor product of the frames of matrix units
in the tensor factors $M_{k^m}(\mathbb{C})$ and $M_{k^{r-m}}(\mathbb{C})$. The matrices from $\pi_1^m(\epsilon)$ commute
with the matrices from $\pi_2^{r-m}(\epsilon)$, moreover, all possible pairwise products of the matrices
from $\pi_1^m(\epsilon)$ by the matrices from $\pi_2^{r-m}(\epsilon)$
(we have exactly $k^{2m}\cdot k^{2(r-m)}=k^{2r}$ such products) give all matrices from the frame $\epsilon.$
If we order the collection of products in the appropriate way, we get the frame $\epsilon.$
The operation which to a pair consisting of commuting $k^m$ and $k^{r-m}$-frames assigns (according to this rule)
the
$k^r$-frame we will denote by dot $\cdot$. In particular, $\beta = \pi_1^m(\beta)\cdot \pi_2^{r-m}(\beta)$ for any $k^r$-frame $\beta.$

Thereby we have defined the continuous maps $\pi_1^m\colon \Fr_{k^r,\, l^s}\rightarrow \Fr_{k^m,\, k^{r-m}l^s}$ and $\pi_2^{r-m}\colon \Fr_{k^r,\, l^s}\rightarrow \Fr_{k^{r-m},\, k^ml^s}.$
In terms of algebra homomorphisms they correspond to the assignment to a homomorphism
$h\colon M_{k^r}(\mathbb{C})\rightarrow M_{k^rl^s}(\mathbb{C})$ its compositions with homomorphisms
$M_{k^m}(\mathbb{C})\rightarrow M_{k^r}(\mathbb{C}),\; X\mapsto X\otimes E_{k^{r-m}}$ and $M_{k^{r-m}}(\mathbb{C})\rightarrow M_{k^r}(\mathbb{C}),\; X\mapsto E_{k^m}\otimes X$
respectively.

\subsection{Functor $\Fr$}

In this subsection we define a functor $\Fr$ from some monoidal category $\mathcal{C}_{k,\, l}$
to the category of topological spaces with a chosen basepoint.

Let us fix an ordered pair of positive integers $k,\, l,\: (k,\, l)=1,\: k,\, l> 1.$
Define the category $\mathcal{C}_{k,\, l}$ whose objects are pairs of the form
$(M_{k^ml^n}(\mathbb{C}),\, \alpha )$,
consisting of a matrix algebra $M_{k^ml^n}(\mathbb{C}),\, m,\, n\geq 0$ and a $k^m$-frame $\alpha$ in it.
A morphism $f\colon (M_{k^ml^n}(\mathbb{C}),\, \alpha )\rightarrow (M_{k^rl^s}(\mathbb{C}),\, \beta )$
is a unital $*$-homomorphism of matrix algebras
$f\colon M_{k^ml^n}(\mathbb{C})\rightarrow M_{k^rl^s}(\mathbb{C})$
such that $f_*(\alpha)=\pi_1^m(\beta),$ where by $f_*$ we denote the map induced on frames by $f$.
Equivalently, we have the commutative diagram
$$
\diagram
M_{k^ml^n}(\mathbb{C}) \rto^f & M_{k^rl^s}(\mathbb{C}) \\
M_{k^m}(\mathbb{C}) \rto^{i_{r,\, m}} \uto^{h_\alpha} & M_{k^r}(\mathbb{C}), \uto_{h_\beta} \\
\enddiagram
$$
where $i_{r,\, m}\colon M_{k^m}(\mathbb{C}) \rightarrow M_{k^r}(\mathbb{C})$ is the homomorphism $X\mapsto X\otimes E_{k^{r-m}}$.

Note that $\mathcal{C}_{k,\, l}$ is a symmetric monoidal category with respect to the bifunctor $\otimes$:
$$
((M_{k^ml^n}(\mathbb{C}),\, \alpha ),\, (M_{k^rl^s}(\mathbb{C}),\, \beta ))\mapsto
(M_{k^ml^n}(\mathbb{C})\otimes M_{k^rl^s}(\mathbb{C}),\, \alpha \otimes \beta)
$$
and the unit object $e:=(M_1(\mathbb{C})=\mathbb{C},\, \varepsilon),$
where $\varepsilon =1$ is the unique $k^0=1$-frame.

Now let us define a functor $\Fr$ from $\mathcal{C}_{k,\, l}$ to the category of topological spaces with a chosen basepoint.
On objects $\Fr(M_{k^ml^n}(\mathbb{C}),\, \alpha )$ is the space of $k^m$-frames in
$M_{k^ml^n}(\mathbb{C})$, where $\alpha$ gives the basepoint. For a morphism
$f\colon (M_{k^ml^n}(\mathbb{C}),\, \alpha )\rightarrow (M_{k^rl^s}(\mathbb{C}),\, \beta )$ put
$$
\Fr(f)\colon \Fr(M_{k^ml^n}(\mathbb{C}),\, \alpha )\rightarrow \Fr(M_{k^rl^s}(\mathbb{C}),\, \beta ),\quad \Fr(f)(\alpha^\prime)=f_*(\alpha^\prime)\cdot \pi_2^{r-m}(\beta).
$$
Then $\Fr(f)$ is a well-defined continuous map preserving basepoints.

Consider a few particular cases.

\smallskip

{\noindent 1)} Suppose $m=0,$ then $\Fr(M_{l^n}(\mathbb{C}),\, \varepsilon )=\{ \varepsilon \}$ is the space consisting of one point, and for a morphism
$f\colon (M_{l^n}(\mathbb{C}),\, \varepsilon )\rightarrow (M_{k^rl^s}(\mathbb{C}),\, \beta )$
the induced map
$$
\Fr(f)\colon \Fr(M_{l^n}(\mathbb{C}),\, \varepsilon )\rightarrow \Fr(M_{k^rl^s}(\mathbb{C}),\, \beta ),\quad \varepsilon \mapsto \varepsilon \cdot \beta =\beta
$$
is the inclusion of the basepoint (note that $\pi^0_1(\beta)=\varepsilon,$ cf. Definition \ref{defframe}, (ii)).

\smallskip

{\noindent 2)} Suppose $n=0,$ then $\Fr(M_{k^m}(\mathbb{C}),\, \alpha )=\PU(k^m)$ and $\alpha$ corresponds to the unit in the group $\PU(k^m)$.
For a morphism $f\colon (M_{k^m}(\mathbb{C}),\, \alpha )\rightarrow (M_{k^r}(\mathbb{C}),\, \beta )$ the diagram
$$
\diagram
\Fr(M_{k^m}(\mathbb{C}),\, \alpha )\rto^{\Fr(f)} \dto_= & \Fr(M_{k^r}(\mathbb{C}),\, \beta ) \dto^= \\
\PU(k^m) \rto^{\ldots \otimes E_{k^{r-m}}} & \PU(k^r) \\
\enddiagram
$$
is commutative (the lower row corresponds to the homomorphism $X\mapsto X\otimes E_{k^{r-m}}$).

\smallskip

{\noindent 3)} $r=m,\; f\colon (M_{k^ml^n}(\mathbb{C}),\, \alpha )\rightarrow (M_{k^ml^s}(\mathbb{C}),\, \beta )$.
Then $\beta =f_*(\alpha),\, \pi_2^{0}(\beta)=\varepsilon$ (cf. Definition \ref{defframe}, (ii)) $\Fr(f)(\alpha^\prime)=f_*(\alpha^\prime).$

\subsection{Natural transformation $\mu \colon \Fr(\ldots ) \times \Fr(\ldots ) \rightarrow \Fr ((\ldots )\otimes (\ldots ))$}

Using the bifunctor $\otimes$ on the category $\mathcal{C}_{k,\, l}$ we define a natural transformation of functors
$\mu \colon \Fr(\ldots ) \times \Fr(\ldots ) \rightarrow \Fr ((\ldots )\otimes (\ldots ))$ from the category
$\mathcal{C}_{k,\, l}\times \mathcal{C}_{k,\, l}$ to the category of topological spaces with a chosen basepoint.
More precisely,
$$
\mu \colon \Fr(M_{k^ml^n}(\mathbb{C}),\, \alpha )\times \Fr(M_{k^pl^q}(\mathbb{C}),\, \varphi )\rightarrow \Fr(M_{k^{m}l^{n}}(\mathbb{C})\otimes M_{k^pl^q}(\mathbb{C})),\, \alpha \otimes \varphi),
$$
$$
\mu(\alpha^\prime ,\, \varphi^\prime)=\alpha^\prime \otimes \varphi^\prime
$$
(recall that $M_{k^{m}l^{n}}(\mathbb{C})\otimes M_{k^pl^q}(\mathbb{C})\cong M_{k^{m+p}l^{n+q}}(\mathbb{C})$),
where $\alpha^\prime \otimes \varphi^\prime$ is the $k^{m+p}$-frame which is the tensor product of the $k^m$-frame $\alpha^\prime$ and the $k^p$-frame $\beta^\prime$.

In fact, $\mu$ is a natural transformation, because for any two morphisms in $\mathcal{C}_{k,\, l}$
$$
f\colon (M_{k^ml^n}(\mathbb{C}),\, \alpha )\rightarrow (M_{k^rl^s}(\mathbb{C}),\, \beta ),\quad
g\colon (M_{k^pl^q}(\mathbb{C}),\, \varphi )\rightarrow (M_{k^tl^u}(\mathbb{C}),\, \psi )
$$
the diagram
\begin{equation}
\label{nattransdiagr}
\diagram
\Fr(M_{k^ml^n}(\mathbb{C}),\, \alpha )\times \Fr(M_{k^pl^q}(\mathbb{C}),\, \varphi )\rto^\mu \dto_{\Fr(f)\times \Fr(g)} &
\Fr(M_{k^{m}l^{n}}(\mathbb{C})\otimes M_{k^pl^q}(\mathbb{C}),\, \alpha \otimes \varphi) \dto^{\Fr(f\otimes g)} \\
\Fr(M_{k^rl^s}(\mathbb{C}),\, \beta )\times \Fr(M_{k^tl^u}(\mathbb{C}),\, \psi )\rto^\mu & \Fr(M_{k^{r}l^{s}}(\mathbb{C})\otimes M_{k^tl^u}(\mathbb{C}),\, \beta \otimes \psi)\\
\enddiagram
\end{equation}
is commutative. Indeed, $\mu \circ (\Fr(f)\times \Fr(g))(\alpha^\prime,\, \varphi^\prime)=\mu(f_*(\alpha^\prime) \cdot \gamma ,\, g_*(\varphi^\prime) \cdot \chi)
=(f_*(\alpha^\prime) \cdot \gamma) \otimes (g_*(\varphi^\prime) \cdot \chi),$
where $\gamma,\, \chi$ are $\pi_2^{r-m}(\beta)$ and $\pi_2^{t-p}(\psi)$ respectively.
On the other hand, $\Fr(f\otimes g)\circ \mu (\alpha^\prime,\, \varphi^\prime)=\Fr(f\otimes g)(\alpha^\prime \otimes \varphi^\prime)=
(f_*(\alpha^\prime) \otimes g_*(\varphi^\prime))\cdot
\Xi,$ where $\Xi$ is the unique $k^{r+t-m-p}$-frame such that
$(f_*(\alpha) \otimes g_*(\varphi))\cdot \Xi=\beta \otimes \psi$.
But $(f_*(\alpha^\prime) \cdot \gamma) \otimes (g_*(\varphi^\prime) \cdot \chi)=(f_*(\alpha^\prime) \otimes g_*(\varphi^\prime))\cdot
(\gamma \otimes \chi)$ by virtue of the commutativity of frames, and moreover
$(f_*(\alpha) \cdot \gamma) \otimes (g_*(\varphi) \cdot \chi)=\beta \otimes \psi.$ Hence
$\Xi =\gamma \otimes \chi$ and the diagram commutes, as claimed.

\subsection{Properties of the natural transformation $\mu$}

First, the natural transformation $\mu$ is associative in the sense that the functor diagram
$$
\diagram
\Fr(\ldots ) \times (\Fr (\ldots ) \times \Fr (\ldots ))\cong (\Fr (\ldots ) \times \Fr (\ldots ))\times \Fr (\ldots ) \dto_{\id \times \mu} \rto^{\quad \qquad \qquad \qquad\mu \times \id} &
\Fr ((\ldots )\otimes (\ldots )) \times \Fr (\ldots ) \dto^\mu \\
\Fr (\ldots )\times \Fr ((\ldots )\otimes (\ldots )) \rto^\mu & \Fr ((\ldots )\otimes (\ldots )\otimes (\ldots ))\\
\enddiagram
$$
commutes (we have used the natural isomorphism $\Fr ((\ldots )\otimes ((\ldots )\otimes (\ldots )))\cong \Fr (((\ldots )\otimes (\ldots ))\otimes (\ldots ))$ in the lower right corner).

Secondly, we need the diagram for identity. Recall that in the monoidal category $\mathcal{C}_{k,\, l}$ $\: e=(\mathbb{C},\, \varepsilon)$ is the unit object, and it is also the initial object.
In particular, for any object $A=(M_{k^ml^n}(\mathbb{C}),\, \alpha )$
there is a unique morphism $\iota_A \colon e\rightarrow A$, i.e. $\iota_A \colon (\mathbb{C},\, \varepsilon)\rightarrow (M_{k^ml^n}(\mathbb{C}),\, \alpha )$.
The identity diagram has the following form:
$$
\diagram
\Fr(e)\times \Fr (\ldots )\rto^{\Fr(\iota) \times \id} \drto & \Fr (\ldots ) \times \Fr (\ldots ) \dto^\mu & \Fr (\ldots )\times \Fr(e) \lto_{\id \times \Fr(\iota)} \dlto \\
& \Fr ((\ldots )\otimes (\ldots )) & \\
\enddiagram
$$
(note that $\Fr(\iota)\colon \Fr(e)\rightarrow \Fr (\ldots )$ is the inclusion of the basepoint). It is easy to see that
(for any pair of objects $A,\, B$ of $\mathcal{C}_{k,\, l}$)
the slanted arrows are homeomorphisms on their images.

There is also the commutativity diagram
$$
\diagram
\Fr (\ldots )\times \Fr (\ldots ) \rrto^\tau \drto_\mu && \Fr (\ldots ) \times \Fr (\ldots ) \dlto^\mu \\
& \Fr ((\ldots ) \otimes (\ldots )) & \\
\enddiagram
$$
(where $\tau$ is the map which switches the factors) which is commutative up to isomorphism. This gives us a homotopy $\mu \circ \tau \simeq \mu$.

For any pair $A,\, B$ of objects of
${\mathcal C}_{k,\, l}$ the natural transformation $\mu$ determines a continuous map $\widetilde{\mu}_{A,\, B}\colon \Fr(A) \times \Fr(B)\to \Fr(A\otimes B)$
of topological spaces.

Put $\Fr_{k^\infty ,\, l^\infty}:=\varinjlim_{\{ f\} }\Fr(M_{k^ml^n}(\mathbb{C}),\, \alpha )$.
Then the above diagrams show that
$\Fr_{k^\infty ,\, l^\infty}$ is a homotopy associative and commutative
$H$-space with multiplication given by $\widetilde{\mu}:=\varinjlim_{A,\, B}\widetilde{\mu}_{A,\, B}$ and with the homotopy unit
$\widetilde{\eta} :=\varinjlim_A \Fr(\iota_A)\colon *=\Fr (e)\rightarrow \Fr_{k^\infty ,\, l^\infty}$.

\subsection{$H$-space structure on the matrix grassmannian}

Note that the analogous construction can be applied to matrix grassmannians (in place of frame spaces).
Namely, consider the category $\mathcal{D}_{k,\, l}$ whose objects are pairs of the form $(M_{k^ml^n}(\mathbb{C}),\, A)$, consisting of a matrix algebra
$M_{k^ml^n}(\mathbb{C}),\; m,\, n\geq 0$ and a $k^m$-subalgebra $A\subset M_{k^ml^n}(\mathbb{C})$ in it
(recall that a $k$-subalgebra is a unital $*$-subalgebra isomorphic $M_{k}(\mathbb{C})$).
A morphism $f\colon (M_{k^ml^n}(\mathbb{C}),\, A)\rightarrow (M_{k^rl^s}(\mathbb{C}),\, B)$ in $\mathcal{D}_{k,\, l}$ is a unital
$*$-homomorphism of matrix algebras $f\colon M_{k^ml^n}(\mathbb{C})\rightarrow M_{k^rl^s}(\mathbb{C})$ such that $f(A)\subset B.$
We define the $k^{r-m}$-subalgebra $C\subset M_{k^rl^s}(\mathbb{C})$ as the centralizer of the subalgebra $f(A)$ in $B$.

Define the functor $\Gr$ from the category $\mathcal{D}_{k,\, l}$ to the category of topological spaces with a chosen basepoint
as follows. On objects the space $\Gr(M_{k^ml^n}(\mathbb{C}),\, A)$ is the space of all $k^m$-subalgebras in $M_{k^ml^n}(\mathbb{C})$ and
$A$ corresponds to its basepoint. For a morphism $f\colon (M_{k^ml^n}(\mathbb{C}),\, A)\rightarrow (M_{k^rl^s}(\mathbb{C}),\, B)$
as above we put
$$
\Gr(f)\colon \Gr(M_{k^ml^n}(\mathbb{C}),\, A)\rightarrow \Gr(M_{k^rl^s}(\mathbb{C}),\, B),\quad \Gr(f)(A^\prime)=f(A^\prime) \cdot C,
$$
where $C$ is the centralizer of the subalgebra $f(A)$ in $B$ and $f(A^\prime) \cdot C$ denotes the subalgebra in $M_{k^rl^s}(\mathbb{C})$
generated by subalgebras $f(A^\prime)$ and $C$ (clearly, $B=f(A)\cdot C$).
Then one can define the analog $\mu^\prime \colon \Gr(\ldots ) \times \Gr(\ldots ) \rightarrow \Gr ((\ldots )\otimes (\ldots ))$ of the natural transformation $\mu$, etc.
(For example, the commutativity of the analog of diagram (\ref{nattransdiagr})
follows from the coincidence $Z_{B\otimes \Psi}(f(A)\otimes g(\Phi))=Z_B(f(A))\otimes Z_\Psi(g(\Phi))$
for any two morphisms $f\colon (M_{k^ml^n}(\mathbb{C}),\, A)\rightarrow (M_{k^rl^s}(\mathbb{C}),\, B)$ and
$g\colon (M_{k^pl^q}(\mathbb{C}),\, \Phi)\rightarrow (M_{k^tl^u}(\mathbb{C}),\, \Psi)$
(which is an analog of the above formula $\Xi =\gamma \otimes \chi$ for frames),
where $Z_B(A)$ denotes the centralizer of a subalgebra $A$ in an algebra $B$.)
This allows us to equip the direct limit $\Gr_{k^\infty,\, l^\infty}:=\varinjlim_{ \{ f\} }\Gr(M_{k^ml^n}(\mathbb{C}),\, A)$ with
the structure of a homotopy associative and commutative $H$-space with a homotopy unit.

Note that there is the functor
$$
\lambda \colon \mathcal{C}_{k,\, l}\rightarrow \mathcal{D}_{k,\, l},\quad (M_{k^ml^n}(\mathbb{C}),\, \alpha)\mapsto
(M_{k^ml^n}(\mathbb{C}),\, M(\alpha)),
$$
where $M(\alpha)$ is the $k^m$-subalgebra spanned on the $k^m$-frame $\alpha$. There is the obvious natural transformation
of functors $\theta \colon \Fr \rightarrow \Gr \circ \lambda$ from the category $\mathcal{C}_{k,\, l}$ to the category
of topological spaces with a chosen basepoint which gives rise to the $H$-space homomorphism
$\Fr_{k^\infty,\, l^\infty}\rightarrow \Gr_{k^\infty,\, l^\infty}$.
Recall that $\Gr_{k^\infty,\, l^\infty}=\Gr \cong \BSU_\otimes$, and the image of the just constructed homomorphism
is the $k$-torsion subgroup in it, as the next proposition claims.

\begin{proposition}
\label{ktors3}
Let $X$ be a compact space. Then the image of the homomorphism
$[X,\, \Fr_{k^\infty,\, l^\infty}]\rightarrow [X,\, \Gr_{k^\infty,\, l^\infty}]$ is the
$k$-torsion subgroup in the group $bsu^0_\otimes(X).$
\end{proposition}
{\noindent \it Proof.\: } This proposition follows from Proposition \ref{ktors2}. Another way is to pass to the direct limit in fibration
$$
\Fr_{k^n,\, l^n}\to \Gr_{k^n,\, l^n}\to \BPU(k^n),
$$
and to notice that the limit map
$\Gr_{k^\infty,\, l^\infty} \to \BPU(k^\infty):=\varinjlim_n\BPU(k^n)$
actually is a localization on $k. \quad \square$

\end{document}